\documentclass[12pt]{amsart}
\usepackage[T1]{fontenc}
\usepackage{amssymb}
\newcommand{\dom}{\operatorname{dom}}

\newcommand{\Aaa }{\mathcal A}
\newcommand{\Raa }{\mathcal R}
\newcommand{\Bee }{\mathcal B}
\newcommand{\Cee }{\mathcal C}
\newcommand{\Wee }{\mathcal W}
\newcommand{\Qee }{\mathcal Q}
\newcommand{\Pee }{\mathcal P}
\newcommand{\Fee }{\mathcal F}

\newcommand{\cl}{\operatorname{cl}}
\renewcommand{\int}{\operatorname{Int}}

\newtheorem{theorem}{Theorem}
\newtheorem{corollary}[theorem]{Corollary}
\newtheorem{lemma}[theorem]{Lemma}

\newtheorem{example}[theorem]{Example}
\newtheorem{remark}[theorem]{Remark}

\author{Andrzej Kucharski}
\address{Andrzej Kucharski \\
 Institute of Mathematics, University of
Silesia \\
 ul. Bankowa 14, 40-007 Katowice}
\email{akuchar@ux2.math.us.edu.pl}
\author{Szymon Plewik}
\address{Szymon Plewik\\Institute of Mathematics,
University of Silesia, ul. Ban\-ko\-wa 14, 40-007 Katowice}
\email{plewik@ux2.math.us.edu.pl}

\begin{document}

\title{Game approach to universally Kuratowski-Ulam spaces}

\begin{abstract}
We consider a version of the open-open game, indicating its connections with universally Kuratowski-Ulam spaces. From \cite{dkz} and \cite{fnr} topological arguments are extracted  to  show that: \textit{Every I-favorable space is universally  Kuratowski-Ulam}, Theorem \ref{FFF};  \textit{If a compact space $Y$ is I-favorable, then the hyperspace $\exp(Y)$ with the Vietoris topology is I-favorable, and hence universally  Kuratowski-Ulam}, Theorems \ref{DDD} and \ref{ooo}.
Notions of uK-U and uK-U$^*$ spaces are compared.
\end{abstract}

\maketitle
\section{Introduction}
     The following theorem was proved (in fact) by K. Kuratowski and S. Ulam, see \cite{ku} and compare \cite{kur} p. 246: \\ \indent 
\textit{Let $X$ and $Y$ be topological spaces such that $Y$ has countable $\pi$-weight. If $E\subseteq X\times Y$ is a nowhere dense set, then there is   $P\subseteq X$ of  first category such that the section $E_x= \{ y: (x,y) \in E \}$
 is nowhere dense in $Y$ for any point $x\in X\setminus P$}.
 
   In   \cite{ox} one can find less general formulation of the Kuratowski Ulam Theorem:\\ \indent 
 \textit{If $E$ is a   plane set of first category, then  $E_x$
 is a linear set of first category for all $x$ except a set of first category.}
 
 In the literature a set of the first category is usually called a meager set. The Kuratowski Ulam  Theorem holds for any meager (nowhere dense) set   $E \subseteq X\times Y$, where the Cartesian product $X\times Y$ is equipped with the  Tychonov topology and   $\pi$-weight of $Y$ is less than additivity of meager sets in $X$, compare \cite{fnr}, \cite{kur} or \cite{ox}. 
 
 The above formulations of the Kuratowski Ulam Theorem suggests two notions of universally Kuratowski-Ulam spaces which one could study.
 
 A space $Y$ is \textit{universally Kuratowski-Ulam} (for short, \textit{uK-U space}), whenever for any topological space $X$ and a meager set $E\subseteq X\times Y$, the set $$ \{x\in X: \{y\in Y: (x,y)\in E\} \mbox{ is not meager in } Y\}$$ is meager in $X$, 
   see D. Fremlin, T. Natkaniec and I. Rec{\l}aw \cite{fnr}.  The class of uK-U spaces has been  investigated in  \cite{fnr}, \cite{fem} and \cite{zsi}.    
 
 A space $Y$ is \textit{universally Kuratowski-Ulam$^*$} (for short, \textit{uK-U$^*$ space}), whenever for a topological space $X$ and a nowhere dense  set $E\subseteq X\times Y$ the set $$ \{x\in X: \{y\in Y: (x,y)\in E\} \mbox{ is not nowhere dense  in } Y\}$$ is meager in $X$, 
   see D. Fremlin \cite{fem}. 
   
   Any uK-U$^*$ space is uK-U space. A proof of this is standard. Indeed, suppose that a space $Y$ is  uK-U$^*$, and $X$  is a topological space. If $E\subseteq X \times Y$ is a meager set, then there exist nowhere dense sets $E^n \subseteq X\times Y$ such that  $E^0 \cup E^1 \cup \ldots \supseteq E$. Put $$P_n = \{x\in X: \{y\in Y: (x,y)\in E^n\} \mbox{ is not nowhere dense  in } Y\}.$$ Each set $P_n$ is meager, hence $P = P_0 \cup P_2 \cup \ldots \subseteq X$ is meager.
    Since  $E_x \subseteq E^0_x \cup E^1_x \cup \ldots $ (recall that $E^n_x= \{y\in Y: (x,y)\in E^n\}$), then $E_x$ is meager  for each  $x\in X\setminus P$.  
   
   The converse is not true: \textit{There is a dense in itself and countable Hausdorff space which is not \mbox{uK-U$^*$}}; see  "6. Examples (b)" in \cite{fem}. 
      Any countable and dense in itself space is meager in itself, and hence has to be uK-U.  The space $C[\omega^\omega]$ of all compact non-empty subsets of the irrationals equipped with the Pixley-Roy topology  has a separable compactification, see A. Szyma\'{n}ski \cite{szy}. One can check that $\omega^\omega\times C[\omega^\omega]$ does not satisfy the Kuratowski Ulam Theorem, hence $C[\omega^\omega]$ is not uK-U$^*$, and any dense subspace of a compactification of $C[\omega^\omega]$ is not \mbox{uK-U$^*$}, too. So, some compactification of $C[\omega^\omega]$  contains a countable Hausdorff space which is  \mbox{uK-U} and not uK-U$^*$. Natural examples of countable
spaces which are not uK-U$^*$ are  spaces of type $Seq$, compare  \cite{vau}.  They are  not uK-U$^*$ by  similar arguments which work with $C[\omega^\omega]$, or with Example 1  in \cite{fnr}.

 The open-open game and I-favorable spaces were  introduced by P. Daniels, K. Kunen and H. Zhou \cite{dkz}.  A space is I-favorable if, and only if it has a club filter, see \cite{dkz}.  Topics of almost the same kind like I-favorable spaces were considered  by  E. V. Shchepin \cite{ksc}, L. Heindorf and L. Shapiro \cite{hs}, and  by B. Balcar, T. Jech and J. Zapletal \cite{bjz}. In \cite{ksc} were introduced $\kappa$-metrizable spaces; in \cite{hs} were considered regularly filtered algebras; in \cite{bjz} were considered semi-Cohen algebras. A Boolean algebra $\mathbb B$ is semi-Cohen (regularly filtered) if, and only if $[\mathbb B]^\omega$ has a closed unbounded set of countable regular subalgebras (contains a club filter). Semi-Cohen algebras and I-favorable spaces are  corresponding classes, compare \cite{bjz}  and \cite{hs}. 
 
 Every dyadic space is uK-U space, see \cite{fnr}. We  extend this fact by showing that any I-favorable space is uK-U$^*$,  Theorem \ref{FFF}. Additionally,   we show  that any hyperspace  $\exp (D^\lambda)$ is uK-U$^*$ space,  Corollary \ref{HHH}.  
\section{The game}
The following game was invented by P. Daniels, K. Kunen  and H. Zhou \cite{dkz}. Two players take turns playing with a topological space $X$. A round consists of Player I choosing a non-empty open set $U\subseteq X$; and Player II choosing a non-empty open set $V\subseteq U$. Player  I wins if the union of all open sets which have been chosen by Player II is dense in $X$. This game was called the \textit{open-open game}.    If the open-open game of uncountable length  is being played  with a space of countable 
cellularity (for example, some $Seq$ spaces), then Player II could be forced to choose disjoint sets at each round. In consequence, Player I wins any such game.  Thus,  any open-open game is not trivial under some restrictions which imply that  Player I can not win always. For example,  rounds are played for each ordinal less than some given ordinal $\alpha$.
  From here,  we consider  cases when  games have the least infinite length i.e. $\alpha = \omega$. 
  
  Let us consider the following game. 
 Player I  chooses  a finite family $\mathcal{A}_0$ of non-empty  open  subsets of $X$. Then Player II chooses a finite family $\mathcal{B}_0$ of non-empty open  subsets of $X$ such that for each $U \in \mathcal{A}_0$ there exists $V \in \Bee_0$ with $V \subseteq U $. Similarly at the $n$-th round  Player I  chooses 
a finite family $\Aaa_n$ of non-empty open  subset of $X$.
Then  Player II chooses a finite family $\mathcal{B}_n$ of non-empty open  subsets of $X$ such that for each $U \in \mathcal{A}_n$ there exists $V \in \Bee_n$ with $V \subseteq U $.
 If for any natural number $k$ the union $\bigcup \{\mathcal{B}_k \cup \mathcal{B}_{k+1} \cup  \ldots\}$ is a dense subset of $X$, then  Player I wins; otherwise Player II wins.

 The space $X$ is I-\textit{favorable} whenever Player I can be insured, by choosing his families $\Aaa_n$ judiciously, that he  wins no matter how Player II plays. In this case we say that  Player I has a winning strategy.  
Player I has a winning strategy whenever  any finite family  of open and disjoint subsets of $X$ he can consider as $\Bee_n$, and then  Player I knows his  $(n+1)$-th round, i.e. he knows how to define $\Aaa_0=\sigma(\emptyset)$ and $\Aaa_{n+1}= \sigma( \Bee_0,\Bee_1 , \ldots, \Bee_n)$.   Any winning strategy  would be defined as  function 
 $$( \Bee_0,\Bee_1 , \ldots, \Bee_n) \mapsto \sigma ( \Bee_0,\Bee_1 , \ldots, \Bee_n), $$ where all families $\Bee_n$ and $\sigma ( \Bee_0,\Bee_1 , \ldots, \Bee_n)$ are finite and  consists of non-empty open sets; and  for any game with succeeding rounds  $ \sigma (\emptyset)$, $\Bee_0$, $\sigma(\Bee_0)$, $\Bee_1,\sigma(\Bee_0, \Bee_1),\ldots, \Bee_n ,\sigma(( \Bee_0,\Bee_1 , \ldots, \Bee_n) $ each union $\bigcup \{\mathcal{B}_k \cup \mathcal{B}_{k+1} \cup  \ldots\}$ is a dense subset of $X$.

 Our definition of I-favorable space is equivalent to the similar definition  stated in  \cite[p. 209]{dkz}.  In fact, if $\mathcal{A}_n = \{ U_1, U_2, \ldots , U_k \}$, then  Player I should play $k$-rounds  choosing  $U_1, U_2, \ldots , U_k$, successively. If Player I has a strategy $\sigma$ which forced Player II to  choose families  $\mathcal{B}_k$ such that  $\bigcup \{\mathcal{B}_0 \cup \mathcal{B}_{1} \cup  \ldots\}$ is a dense subset of $X$, then Player I could divide the set of natural numbers onto infinite many of pairwise disjoint infinite pieces. Then Player I could play at each piece following $\sigma$, and he obtains the winning strategy.   In consequence,  for the definition of  I-favorable spaces one can use  the open-open game, or the topological version of the game $G_4$, see \cite[p. 219]{dkz}.

Many cases when Player II can be insured that he  wins no matter how Player I plays were considered in \cite{dkz} or \cite{szy}. By Theorem \ref{FFF} spaces $Seq$ are not I-favorable. However, one can check directly that  Player II could always win a game with any $Seq$: Any $Seq$ has a tiny sequence, compare \cite{szy}, and therefore Player II has winning strategy.   

Let us recall a few comments according to \cite{dkz}. 
Any space with countable $\pi$-weight is I-favorable. Indeed, if $\{W_0, W_1, \ldots \}$ is a  
$\pi$-base for $X$, then Player I chooses $\mathcal{A}_n$ such that always there exists $U \in \mathcal{A}_n$ and $U \subseteq  W_n$.
  If a space $X$ has uncountable cellularity, then $X$ is not \mbox{I-favorable.} Indeed, there exists an uncountable family $\mathcal{W}$ of open and disjoint subsets of $X$,  and Player II can choose $\mathcal{B}_n$ such that always  $\bigcup \mathcal{B}_n$ intersects  finite many members of  $\mathcal{W}$. Another example is a regular Baire space $ X$ with a category measure $\mu$ such that $\mu (X) =1$ (for more details  see \cite[p. 86 - 91]{ox}). Any such  $X$ is not I-favorable, since Player II can choose $\mathcal{B}_n$ such that always  $\mu (\bigcup \mathcal{B}_n) < \frac{1}{2^{n+2}}$. This follows  $\mu \left(X \setminus (\bigcup \{\mathcal{B}_0 \cup \mathcal{B}_1 \cup  \ldots\})\right) \geq \frac{1}{2}$. Therefore  the complement $X \setminus (\bigcup \{\mathcal{B}_0 \cup \mathcal{B}_1 \cup  \ldots\})$ has to have non empty interior.
\section{On I-favorable spaces}
A topological characterization of I-favorable spaces  is applied to describe direct proofs of some know facts. Moreover, we show that if a compact space $X$ is I-favorable, then the hyperspace $\exp(X)$ with the Vietoris topology is I-favorable.
 We extract topological versions of arguments used  in \cite{dkz} and \cite{fnr}.

For any Cantor cube $D^\lambda$ fix the following notation. Let $\lambda$ be a cardinal number, $D = \{0,1\}$, and let $D^\lambda$ be equipped with  the product topology. The product topology is generated by subsets  $ \{ q \in D^\lambda: q(\alpha) = k\},$ where $\alpha \in \lambda$ and $k\in D$. If $f: Y \to D$ and $Y\in [\lambda]^{<\omega}$,  then  $W_f = \{ q \in D^\lambda:  f \subseteq q\} $. All sets $W_f$ constitute  an open base.
\begin{example} \label{can} The Cantor cube $D^{\lambda}$ is I-favorable.
\end{example}
\begin{proof}
Player I put $\mathcal{A}_0 = \{ D^\lambda\}$. If a family $\Bee_0$  is defined, then Player I chooses  base  open sets 
$W_q \subseteq Q$ for any $Q\in \Bee_0$ and put \mbox{$\Aaa_1 = \{ W_f: f\in D^{J_1}\}$,} where  $J_1= \bigcup \{\dom{(q)}: W_q \subseteq Q \in \Bee_0 \}.$ Player I wins, whenever at the $n$-th round he always chooses  base  sets 
$W_q \subseteq Q$ for any $Q\in \Bee_{n-1} $, and  put $\Aaa_n = \{ W_f: f\in D^{J_n}\}$, where $J_n= \bigcup \{\dom{(q)}: W_q \subseteq Q \in \Bee_{n-1} \}.$
Any such played game defined a sequence $J_1 \subseteq J_2 \subseteq \ldots $ of finite subsets of $\lambda$. 
Fix a base set $W_f$ where $f \in D^J$, i.e. $J = \dom(f)$. Take a natural number $n$ such that $J\cap J_n = J \cap J_{n+1}$, and next take
$q\in D^{J_n}$ such that  functions $f$ and $q$ are compatible on the set $J\cap J_n = \dom (f) \cap \dom (q)$. There exists $q^* \in D^{J_{n+1}}$ such that  
$$ \Aaa_n \ni W_q \supseteq V \supseteq W_{q^*} \in \Aaa_{n+1},$$ where $V \in \Bee_n$. Functions $f$ and $q^*$ are compatible on the set $$J\cap J_n = \dom (f) \cap \dom (q^*)= J\cap J_{n+1}.$$ Therefore $W_f$ meets $W_{q^*}$, and hence $\emptyset \neq W_f \cap W_{q^*} \subseteq W_f \cap V \subseteq V$. Since $n$    could be arbitrarily large and  $V \in \Bee_n$, then  each  $\bigcup \{\mathcal{B}_k \cup \mathcal{B}_{k+1} \cup  \ldots\}$ has to be a dense subset of $X$.
\end{proof}
 We have repeated a  special case of Theorem 1.11, see  \cite{dkz}. Our proof of Example \ref{can} explicitly defines  a winning strategy.
But if families $\Aaa_0, \Aaa_1, \ldots$  have been defined simultaneously, then Player I would lose. This would not happen when $X$ has countable $\pi$-base. However  for  \mbox{$X=D^\lambda$,} where $\lambda $ is uncountable, this is possible. Indeed, if  Player I fixes each family $\Aaa_n$, then Player II could choose a finite family $\Bee^*_n$ such that for any $U\in \Aaa_n$ there exists a base subset $W_q \in \Bee^*_n$ with $W_q\subseteq U$.  Put $J_n= \bigcup \{\dom{(q)}: W_q \in  \Bee^*_n\}$, and  take an index \mbox{$\alpha \in \lambda \setminus(J_0 \cup J_1 \cup \ldots )$.}  Afterwards Player II put  $$\Bee_n =  \{V \cap \{q\in D^\lambda: q(\alpha) =1\}:  V \in \Bee^*_n \}.$$  No member of $\Bee_n$ meets  $\{q\in D^\lambda: q(\alpha) =0\}$. In fact, we get the following.
\begin{remark}\label{rem} For each sequence $(\Aaa_0, \Aaa_1, \ldots)$ consisting of finite non-empty families of open subsets of $D^\lambda$, there is a corresponding sequence $(\Bee_0, \Bee_1, \ldots)$ consisting of finite non-empty families $\Bee_n$ of non-empty open sets such that each $\Bee_n$ refines $\Aaa_n$, and yet the union $\Bee_0 \cup \Bee_1 \cup \ldots$ is not dense. \qed
\end{remark}
 Countable subsets of $\lambda$ are important in our proof of Example~ \ref{can}. Any $J\in[\lambda]^{\omega}$ fixes the countable family of base sets $$ C_J=\{W_f: f:Y\to D \mbox{ and } Y\in[J]^{<\omega}\},$$
which  fulfills the following condition: 
 
\textit{For any open  $V \subseteq D^{\lambda}$ there is $W \in C_J$ such that if $U \in C_J$ and $U \subseteq W$, then $U \cap V \not= \emptyset $}.

 \noindent This condition may be considered in an arbitrary topological space $X$ with a fixed $\pi$-base $\Qee$.
 According to definitions \cite[p. 208]{dkz} a family $\Cee \subset [\Qee]^{\omega}$ is called \textit{a club filter} whenever:

(1) \textit{The family $\Cee$ is closed under $\omega$-chains with respect to inclusion, i.e. if $\Pee_1 \subseteq \Pee_2 \subseteq  \ldots $ is an  $\omega$-chain which consists of elements of $\Cee$, then  $\Pee_1 \cup \Pee_2 \cup \ldots \in \Cee$};

(2) \textit{For any countable subfamily  $\Aaa\subseteq \Qee$, where $\Qee$ is the   $\pi-$base fixed above,  there  exists $\Pee \in \Cee$ such that   $ \Aaa\subseteq\Pee$};

(3)  \textit{For any non-empty open set $V$ and each $\Pee \in \Cee$  there is $W \in \Pee$ such that if $U \in \Pee$ and $U \subseteq W$, then $U$ meets $V$, i.e. $U \cap V \not= \emptyset $}. 

  Conditions $(1) - (3)$ are extracted from properties of Cantor cubes used in Example \ref{can}. The following two lemmas repeat Theorem 1.6, see \cite{dkz}. 
\begin{lemma} \label{BBB}  If a topological space has a club filter, then it is I-favorable.
\end{lemma}
\begin{proof}
 Without lost of generality one can assume that any $\Bee_n$ will be contained in $\Qee$. Let $\Aaa_0 =\{ X\}$.
 If $\Bee_0$ has been defined, then  Player  I chooses $\Pee_0\in\Cee$ such that  $\Bee_0 \subseteq \Pee_0$, by (2).
Enumerate  $\Pee_0=\{V_0^0,V_1^0,\ldots \}$ and put $\Aaa_1=\{V_0^0\}.$
If families $\Bee_n$ and $\Pee_{n-1}$ have been defined, then Player I chooses $\Pee_n\in\Cee$ such that $\Bee_n\cup \Pee_{n-1}\subseteq \Pee_{n}$, using (2) again.
Let $\Pee_n=\{V_0^n,V_1^n,\ldots \}$, and put $\Aaa_{n+1}=\{V_j^i: i\leqslant n \mbox{ and } j\leqslant n \}$. By Condition $(1)$, let $\Pee_0\cup\Pee_1\cup\ldots =\Pee_\infty \in \Cee$. 
We shall show that any union $\bigcup \{\mathcal{B}_k \cup \mathcal{B}_{k+1} \cup  \ldots\}$ is a dense subset of $X$.  Suppose that $V$ is a non-empty open set  such that $  V \cap\bigcup \{\mathcal{B}_k \cup \mathcal{B}_{k+1} \cup  \ldots\} =\emptyset . $ By (3) choose $V_j^i\in \Pee_\infty$ such that if $U \in \Pee_\infty$ and $U \subseteq V_j^i$, then $U \cap V \not= \emptyset $.
Take  $m\geqslant\max\{i,j,k\}$. There exists $W\in\Bee_{m+1} \subseteq \Pee_\infty$ such that $W\subseteq V_j^i$, hence $W\cap V\ne\emptyset$. But $W\in  \mathcal{B}_k \cup \mathcal{B}_{k+1} \cup  \ldots$, a contradiction. 
\end{proof}
\begin{lemma} \label{AAA} If a topological space is  I-favorable, then it has a club filter such that any of its elements is closed under finite intersection.
\end{lemma}
\begin{proof}
Let $\Qee$ be a fixed $\pi$-base,  which is closed under finite intersection, and let $\sigma $ be a winning strategy for Player I. For each countable family $\Raa\in[\Qee]^{\leqslant\omega}$ let $\Raa_1$ be the closure under finite intersection of $\Raa$ and the family $$ \bigcup \{\sigma( \Fee_0,\Fee_1 , \ldots, \Fee_k):\{ \Fee_0,\Fee_1 , \ldots, \Fee_k\} \subset [\Raa]^{<\omega}\mbox{ and } k \in \omega \}.$$   By induction, let $\Raa_{n+1}$ be the closure under finite intersection of $\Raa_{n}$ and $$  \bigcup \{\sigma( \Fee_0,\Fee_1 , \ldots, \Fee_k):\{ \Fee_0,\Fee_1 , \ldots, \Fee_k\} \subset[\Raa_n]^{<\omega}\mbox{ and } k \in \omega \}.$$   A desired club filter $\Cee$ consists of all unions $\Raa_1\cup\Raa_2\cup\ldots$, where $\Raa\in[\Qee]^{\leqslant\omega}$.
By the definition any element of $\Cee$ is closed under finite intersection. Consider an  $\omega$-chain $\Pee_1 \subseteq \Pee_2 \subseteq  \ldots $   in $\Cee$. Let \mbox{$   \Pee_1 \cup \Pee_2 \cup   \ldots =\Raa$.} If $\Fee_0,\Fee_1 , \ldots, \Fee_k$ are    finite families contained in $\Raa$, then there exists  $n$ such that $\Fee_0 \cup \Fee_1 \cup \ldots \cup \Fee_k \subseteq \Pee_n $ and $\sigma(\Fee_0,\Fee_1 , \ldots, \Fee_k)\subseteq \Pee_{n+1}$. This follows $ \Raa \in \Cee$, i.e. Condition (1) holds.  Condition (2) follows directly from  the definition of $\Cee$. Suppose that $\Pee \in \Cee$ and an open set $V$ fulfill the negation of (3). Then, Player II chooses families consisting of sets disjoint with $V$. In consequence, he wins the game $\sigma(\emptyset), \Bee_0, \sigma(\Bee_0), \Bee_1, \ldots$, a contradiction.
\end{proof}
The next corollary was proved in  \cite[Corollary 1.7]{dkz}. 
\begin{corollary} Any product of I-favorable spaces is I-favorable.
\end{corollary}
\begin{proof}
Consider a product $\prod\{X_\alpha:\alpha\in T\}$, where any $X_\alpha$
is I-favorable. Let $\Cee_\alpha$ be a club filter which witnesses that $X_\alpha$ is I-favorable, where  $\Qee_\alpha$ is a $\pi$-base needed in Condition (2). Fix $\lambda\in [T]^{\omega}$ and $\Pee_\alpha\in
\Cee_\alpha$ for each $\alpha\in\lambda$. Let $\Pee(\lambda)$ be the family of all  $\prod\{W_\alpha:\alpha\in S\}$, where $S\in[\lambda]^{<\omega}$ and $W_\alpha\in\Pee_\alpha\in\Cee_\alpha$.
The family $\Cee=\{\Pee(\lambda):\lambda\in [T]^{\omega}\}$ is a desired club filter. 
\end{proof}
In \cite[p. 210]{dkz} it was proved that dyadic spaces are I-favorable. But L. Shapiro \cite{s76} show that some hyperspaces over dyadic spaces can be non-dyadic. For example,  $\exp( D^{\omega_2})$ is a non-dyadic space. For some facts and notions concerning a hyperspace with the Vietoris topology, which are not defined here, see \cite{kur}.  Now,  prove  the following. 
\begin{theorem}\label{DDD} If a compact space $X$ is I-favorable, then the hyperspace $\exp(X)$ with the Vietoris topology is I-favorable, too.
\end{theorem}
\begin{proof}
Fix a $\pi$-base $\Qee$ closed under finite intersection,  and a club filter $\Cee$ for $X$. If $n$ is a natural number and $V_1,V_2,\ldots , V_n$ are open subsets of $X$, then let $<V_1,V_2,\ldots , V_n>$ denotes the family of all closed sets $A\subseteq V_1\cup V_2\cup \ldots \cup V_n $ such that $A\cap V_i\ne\emptyset$ for $1\leqslant i\leqslant n$. The family $$\Qee^*=\{<V_1,V_2,\ldots , V_n>:V_i\in \Qee \mbox{ for } 1\leqslant i\leqslant n\}$$ is  a $\pi$-base for $\exp(X)$.  For any $\Pee\in\Cee$,  let $$\Pee^*=\{<V_1,V_2,\ldots , V_n>:V_i\in \Pee \mbox{ for } 1\leqslant i\leqslant n\}.$$ We   shall check that the family $\Cee^*=\{\Pee^*:\Pee\in\Cee\}$ is a club filter for $\exp(X)$. Then the result follows from Lemma \ref{BBB}.  

By definitions $\Cee^*$ fulfills conditions $(1)$ and $(2)$ and any family $\Pee^* \in \Cee^*$ is closed under finite intersection. Consider an open set $<V_1,V_2,\ldots , V_n> \subseteq \exp(X)$ and a family $\Pee \in \Cee$. For  $1\leqslant i \leqslant n$, by   $(3)$, choose $W_i \in \Pee$ such that if $U\in \Pee$ and $U\subseteq W_i$, then $U$ meets $ V_i$. If $$ 
<W_1,W_2,\ldots , W_n> \supseteq <U_1,U_2,\ldots , U_m> \in \Pee^*,$$ then fix $U^j_i \in \{ U_1,U_2,\ldots , U_m\}$  with $U^j_i \subseteq W_i$.  Since $\Pee$ is closed under finite intersection, then  $U^j_i \cap W_i \in \Pee$. By $(3)$  choose $x_i \in V_i \cap W_i \cap U^j_i$ for $1\leqslant i \leqslant n$. Similarly, choose $y^j_i\in V_i \cap U_j \cap W_i$ whenever  $U_j$ meets $W_i$.  The closed (finite)  set $$ \{x_i: 1\leqslant i \leqslant n\} \cup \{y^j_i: 1\leqslant j \leqslant m \mbox{ and } 1\leqslant i \leqslant n \} \subseteq X$$ belongs to the intersection $<V_1,V_2,\ldots , V_n> \cap <U_1,U_2,\ldots , U_m>$. It follows that  (3) holds for  $\exp (X)$. 
\end{proof}
Special cases of Theorem \ref{DDD} could be  deduced in another way. 
L. Shapiro  observed that $\exp(D^{\lambda})$ is co-absolute with $D^\lambda$, see \cite[Theorem 4]{s85} and \cite[p.17-18]{sce}.  Therefore one could obtain that $\exp(D^{\lambda})$ is I-favorable by \cite[Fact 1.3]{dkz}.

One can check that if there is a club filter $\Cee$  for $\exp(X)$ such that any $\Pee \in \Cee$ consists of base sets of the form $<V_1, V_2, \ldots , V_n>$, then families  constitute all  $V_i$ such that $V_i \in \{ V_1, V_2, \ldots , V_n\} $, where $ <V_1, V_2, \ldots , V_n> \in \Pee $ consists of a club filter for $X$. This gives the converse of Theorem \ref{DDD}. 
\section{On uK-U$^*$ spaces}
 In this note  the next theorem is  main novelty.  Closed nowhere dense sets are valid for uK-U$^*$ properties. Now, it will be  convenient for us  to use 
 open and dense subsets of $X\times Y$, instead of nowhere dense ones. In the  proof of Theorem \ref{EEE} Player II uses an obvious fact:  \textit{If a dense subset $E\subseteq X\times Y$ is open, then for any non-empty open sets $U$ of $ X$ and  $V_1, V_2, \ldots V_n$ of $Y$ there exist non-empty open sets $U^*\subseteq U$ and \mbox{$V_1^*\subseteq V_1, V_2^*\subseteq V_2 \ldots V_n^*\subseteq V_n $}  such that always $U^* \times V^*_i \subseteq E$.} 
\begin{theorem}\label{EEE} Suppose $X$ and $Y$ are topological spaces, where $Y$ is I-favorable. If a set $E\subseteq X\times Y$ is open and dense with respect to the product topology, then there exists a meager subset $P\subseteq X$ such that   the section $$E_x= \{ y\in Y: (x,y) \in E \}$$ is dense in $Y$ for  all $x \in X \setminus P$.
\end{theorem}
\begin{proof}
If Player I has chosen  a finite family $\Aaa_0$ of open and disjoint subsets of $Y$, then Player II chooses an open set $Q_0 \subseteq X$ and a finite  family $\Bee_0(Q_0)$ of open and disjoint subset of $Y$ such that for each $U\in \Aaa_0$ there exists $V\in \Bee_0(Q_0)$ with $V\subseteq U$ and $Q_0 \times V \subseteq E$.

 Afterwards  Player I chooses a finite family $\Aaa_1(Q_0)$  of open and disjoint subsets of $Y$  in accordance with to his winning strategy at the round following after $\Aaa_0$, $\Bee_0(Q_0)$. 
  
  Assume that open sets $X \supseteq Q_0 \supseteq Q_1 \supseteq \ldots \supseteq Q_{n-1}$ and finite families $\Aaa_0,  \Bee_0(Q_0),\Aaa_1(Q_0), \ldots,\Bee_{n-1}(Q_{n-1}), \Aaa_n(Q_{n-1})$ are defined. Then Player II chooses an open set $Q_n \subseteq Q_{n-1}$ and a finite  family $\Bee_n(Q_n)$ of open and disjoint subset of $Y$ such that for each $U\in \Aaa_n(Q_{n-1})$ there exists $V\in \Bee_n(Q_n)$ with $V\subseteq U$ and $Q_n \times V \subseteq E$. 
  
  Afterwards  Player I chooses a finite family $\Aaa_{n+1}(Q_n)$  of open and disjoint subsets of $Y$ in accordance with his winning strategy in the round following after $\Aaa_0, \Bee_0(Q_0),   \ldots , \Bee_{n-1}(Q_{n-1}),\Aaa_n(Q_{n-1}), \Bee_n(Q_n)$.

Let $\Wee_0$ be some maximal family of open and disjoint subsets of $X$ from which Player II could  choose at start as  sets  $Q_0$.  Suppose that families $\Wee_0, \Wee_1, \ldots ,\Wee_{n-1}$ are defined. 
Let $\Wee_n^Q$ be a maximal family of open and disjoint subsets of $X$ which Player II could  choose    at the round following after $\Aaa_0,  \Bee_0(Q_0), \ldots,\Bee_{n-1}(Q_{n-1}), \Aaa_n(Q_{n-1})$,  where $ Q_0 \supseteq Q_1 \supseteq \ldots \supseteq Q_{n-1}$ and $Q_i \in \Wee_i$, for $0\leq i \leq n-1$.  Put $$\Wee_n = \bigcup \{\Wee_n^Q: Q \in \Wee_{n-1}\}.$$
By the induction families $\Wee_0, \Wee_1, \ldots$ are defined. 
Any $\bigcup \Wee_n$ is an open dense subset of $X$.
If always $Q_n \in \Wee_n$ and $x\in Q_0 \cap Q_1\cap \ldots$, then any union $$\bigcup \{ \Bee_k(Q_k) \cup \Bee_{k+1}(Q_{k+1}) \cup \ldots \} $$ is a dense subset of $Y$ since the winning strategy of I forces moves   $\Bee_0(Q_0),  \Bee_1(Q_1),  \ldots $ with a such property. But $V\in \Bee_n(Q_n)$ implies $Q_n \times V \subseteq E$. Therefore $E_x$ should be dense in $Y$. Families
$\Wee_n$ are maximal and consists of open sets, so $\bigcup \Wee_n$ is always open and dense in $X$. Hence for any $$x\in \bigcup \Wee_0 \cap \bigcup \Wee_1 \cap \ldots$$ the set $E_x$ should be dense in $Y$.
Let $P=X\setminus (\bigcup \Wee_0 \cap \bigcup \Wee_1 \cap \ldots)$.
\end{proof}

Apply the above theorem to indicate connections between games and  universally Kuratowski-Ulam spaces.
 
 \begin{theorem} \label{FFF}  Every  I-favorable space is uK-U$^*$.
\end{theorem}
\begin{proof} Suppose that a space $Y$ is  I-favorable, and $X$  is a topological space. If $D\subseteq X \times Y$ is nowhere dense, then it's closure is nowhere dense, too. Apply Theorem \ref{EEE} with $E=X\times Y \setminus \cl D$. 
\end{proof}

Thus, there has been  given an argument which suggests that an adequate   meaning of universally Kuratowski-Ulam spaces should be uK-U$^*$ spaces, compare \cite{fem}.   There exist non-dyadic and compact spaces which are  uK-U$^*$. 
 \begin{theorem} \label{ooo} If a compact space $Y$ is I-favorable, then the hyperspace $\exp(Y)$ with the Vietoris topology is uK-U$^*$.
\end{theorem}
\begin{proof}  The hyperspace $\exp (Y)$ is  I-favorable by Theorem \ref{DDD}. So, one could apply Theorem  \ref{FFF}. 
\end{proof}
\begin{corollary} \label{HHH} If $\lambda > \omega_1$, then the hyperspace $\exp(D^\lambda)$ is uK-U$^*$ and non-dyadic.
\end{corollary} 
\begin{proof}
For any cardinal $\lambda > \omega_1$ the hyperspace  $\exp(D^\lambda)$ is non-dyadic, by \cite{s76}. The Cantor cube $D^\lambda$ is I-favorable and hence $\exp (D^\lambda)$ is I-favorable by Theorem \ref{DDD}. Theorem \ref{ooo} implies that $\exp(D^\lambda)$ is uK-U$^*$. 
\end{proof}
\section{Final remarks}
In \cite[Theorem 1]{s85}  
L. Shapiro  showed that any dyadic space  is co-absolute with a finite disjoint union of Cantor cubes or is co-absolute with the one point compactification of countable many Cantor cubes. Therefore, any dyadic space is  co-absolute with some I-favorable space. One can check this using the definition of I-favorable space. So, one  can reprove \cite[Theorem 1.11]{dkz} using \cite[Fact 1.3]{dkz}. In other words, any dyadic space is I-favorable since it is co-absolute with a I-favorable space. This and  Theorem \ref{FFF} give a proof that  dyadic spaces are universally Kuratowski-Ulam. We have reproved Corollary 3 from \cite{fnr}.  Similarly, by  Theorem \ref{FFF}, and Corollary 5.5.5 \cite{hs}, and Proposition 5.5.6 \cite{hs} one obtains that any space which is co-absolute with a $\kappa$-metrizable space is  uK-U$^*$, compare \cite{ksc},  \cite[p. 44]{hs}. However, we do not know: \textit{Does there exist  a compact universally Kuratowski-Ulam space which is not 
 I-favorable}?
 
 \textbf{Acknowledgments}. We  thank  M. Hru\v{s}\'{a}k for  a helpful remark  on relations between  I-favorable spaces and semi-Cohen Boolean algebras. We are indebted to the referee for very careful reading of the  paper. His (or her) constructive comments drastically revise the previous version of this note. Remark \ref{rem} was suggested by the Referee. In the present version  the notion of uK-U$^*$ is considered. Previously, uK-U spaces have been commented, only. The references  \cite{fem}, \cite{ksc}, \cite{sce}, \cite{szy}, \cite{vau} and \cite{zsi} are  enclosed, now.

\end{document}